\documentclass[11pt]{article}

\usepackage{amsmath, amssymb, amsthm, amstext, graphicx,tikz}
\usepackage{enumerate}
\usepackage{fullpage}

\usepackage{color}

\newtheorem{theorem}{Theorem}[section]

\newtheorem{lemma}[theorem]{Lemma}

\theoremstyle{definition}

\def\e{\epsilon}
\def\G{\mathcal{G}}

\title{Combinatorial theorems relative to a random set}
\author{David Conlon\thanks{Mathematical Institute, Oxford OX2 6GG,
United Kingdom. Email: {\tt david.conlon@maths.ox.ac.uk}. Research supported by a Royal Society University Research Fellowship.}}
\date{}

\begin{document}

\maketitle

\begin{abstract}
We describe recent advances in the study of random analogues of combinatorial theorems.
\end{abstract}

\section{Introduction} \label{sec:intro}

Random graphs have played an integral role in extremal combinatorics since they were first used by Erd\H{o}s \cite{E47} to prove an exponential lower bound for Ramsey numbers. The binomial random graph $G_{n,p}$ is a graph on $n$ vertices where each of the $\binom{n}{2}$ possible edges is chosen independently with probability $p$. In modern terminology, Erd\H{o}s' result says that with high probability $G_{n,1/2}$ contains no clique or independent set of order $2 \log_2 n$. This then translates to a lower bound of $2^{t/2}$ for the Ramsey number $R(t)$ (this will be defined in Section~\ref{sec:Ramsey}). 

Although there were several applications of random graphs prior to their work, the first systematic study of random graphs was undertaken by Erd\H{o}s and R\'enyi \cite{ER59, ER60}. The concept of random graphs was also introduced independently by several other authors but, as explained by Bollob\'as \cite{B01}, `the other authors were all concerned with enumeration problems and their techniques were essentially deterministic.' Though it has its origins in applications to extremal combinatorics, the theory of random graphs is now a rich and self-sustaining area of study (see, for example, \cite{B01, JLR00}).

Suppose that $\mathcal{P}$ is a graph property, that is, a family of graphs closed under isomorphism. In studying random graphs, we are usually concerned with determining the probability that $G_{n,p}$ is in $\mathcal{P}$ for some property $\mathcal{P}$. For many properties, this probability exhibits a phase transition as $p$ increases, changing abruptly from $0$ to $1$. The crossover point is known as the threshold. Formally, we say that $p^*:= p^*(n)$ is a threshold for $\mathcal{P}$ if 
\[
\lim_{n \rightarrow \infty} \mathbb{P} [G_{n,p} \mbox{ is in $\mathcal{P}$}] =
\begin{cases}
0 & \text{if $p = o(p^*)$}, \\
1 & \text{if $p = \omega(p^*)$}.
\end{cases} 
\]
We note that, depending on the property $\mathcal{P}$, the probability could also collapse from $1$ to $0$ as $p$ increases. However, for most properties considered in this paper, the behaviour is as above. To give some simple examples, the properties of being connected and having a Hamiltonian cycle are both known to have a threshold at $p^* = \frac{\log n}{n}$, while the property of containing a particular graph $H$ has a threshold at $n^{-1/m(H)}$, where
\[m(H) = \max\left\{\frac{e(H')}{v(H')}: H' \subseteq H\right\}.\]
This function reflects the fact that a graph appears once its densest subgraph does.

Since the late eighties, there has been a great deal of interest in determining thresholds for analogues of combinatorial theorems to hold in random graphs and random subsets of other sets such as the integers. To give an example, we say that a graph $G$ is $K_3$-Ramsey if any $2$-colouring of the edges of $G$ contains a monochromatic triangle. One of the foundational results in this area, proved by Frankl and R\"odl \cite{FR86} and \L uczak, Ruci\'nski and Voigt \cite{LRV92},  then states that there exists $C > 0$ such that if $p > C/\sqrt{n}$ then
\[\lim_{n \rightarrow \infty} \mathbb{P} [G_{n,p} \mbox{ is $K_3$-Ramsey}] =1.\]
Frankl and R\"odl used this theorem to prove that there are $K_4$-free graphs which are $K_3$-Ramsey, a result originally due to Folkman \cite{F70}. However, this new method allowed one to prove reasonable bounds for the size of such graphs, something which was not possible with previous methods.

From this beginning, a large number of papers were written on sparse random analogues of combinatorial theorems. These included papers on analogues of Ramsey's theorem, Tur\'an's theorem and Szemer\'edi's theorem, though in many cases these efforts met with only partial success. This situation has changed dramatically in recent years and there are now three distinct, general methods for proving sparse random analogues of combinatorial theorems, furnishing solutions for many of the outstanding problems in the area. 

The first two of these methods were developed by Gowers and the author \cite{CG14} and, independently, by Schacht \cite{S14} and Friedgut, R\"odl and Schacht \cite{FRS10}. The third method was found later by Balogh, Morris and Samotij \cite{BMS14} and, independently, by Saxton and Thomason \cite{ST14}. 
Broadly speaking, the method employed by Gowers and the author builds on the transference principle developed by Green and Tao \cite{GT08} in their proof that the primes contain arbitrarily long arithmetic progressions; the method of Schacht and Friedgut, R\"odl and Schacht extends a multi-round exposure technique used by R\"odl and Ruci\'nski \cite{RR95} in their study of Ramsey's theorem in random graphs; and the third method is a byproduct of general results about the structure of independent sets in hypergraphs, themselves building on methods of Kleitman and Winston \cite{KW82} and Sapozhenko \cite{S87, S01, S06}. Of course, this summary does a disservice to all three methods, each of which involves the introduction of several new ideas. Surprisingly, all three proofs are substantially different and all three methods have their own particular strengths, some of which we will highlight below.

Rather than focusing on these three methods from the outset, we will further describe the developments leading up to them, explaining how these new results fit into the broader context. This will also allow us to review many of the important subsequent developments. We begin by discussing random analogues of Ramsey-type theorems. 

\section{Ramsey-type theorems in random sets} \label{sec:Ramsey}

Ramsey's theorem \cite{R30} states that for any graph $H$ and any natural number $r$ there exists $n$ such that any $r$-colouring of the edges of the complete graph $K_n$ on $n$ vertices contains a monochromatic copy of $H$. The smallest such $n$ is known as the $r$-colour Ramsey number of $H$ and denoted $R(H; r)$. When $r = 2$, we simply write this as $R(H)$ and when $H = K_t$, we just write $R(t)$. The result of Erd\H{o}s mentioned in the introduction then says that $R(t) \geq 2^{t/2}$, while an upper bound due to Erd\H{o}s and Szekeres \cite{ES35} says that $R(t) \leq 4^t$. Though there have been lower order improvements to both of these estimates \cite{C09, S75}, it remains a major open problem to give an exponential improvement to either of them.

Given a graph $H$ and a natural number $r$, we say that a graph $G$ is $(H,r)$-Ramsey if in any $r$-colouring of the edges of $G$ there is guaranteed to be a monochromatic copy of $H$. Ramsey's theorem is the statement that $K_n$ is $(H, r)$-Ramsey for $n$ sufficiently large, while the overall aim of graph Ramsey theory is to decide which graphs are $(H,r)$-Ramsey for a given $H$ and $r$. Though coNP-hard in general \cite{B90}, this problem has borne much fruit and there is now a large theory with many interesting and important results (see, for example, \cite{GRS90}). One of the highlights of this theory is the following random Ramsey theorem of R\"odl and Ruci\'nski \cite{RR93, RR94, RR95}, which determines the threshold for Ramsey's theorem to hold in random graphs. As mentioned in the introduction, this result built on earlier work of Frankl and R\"odl \cite{FR86} and \L uczak, Ruci\'nski and Voigt \cite{LRV92}. Here and throughout the paper, we will write $v(H)$ and $e(H)$ for the number of vertices and edges, respectively, of a graph $H$.

\begin{theorem} \label{thm:RR}
For any graph $H$ that is not a forest consisting of stars and paths of length $3$ and any positive integer~$r \geq 2$, there exist positive constants $c$ and $C$ such that 
\[
\lim_{n \rightarrow \infty} \mathbb{P} [G_{n,p} \mbox{ is $(H,r)$-Ramsey}] =
\begin{cases}
0 & \text{if $p < c n^{-1/m_2(H)}$}, \\
1 & \text{if $p > C n^{-1/m_2(H)}$},
\end{cases}
\]
where 
\[m_2(H) = \max\left\{\frac{e(H') - 1}{v(H') - 2}: H' \subseteq H \mbox{ and } v(H') \geq 3\right\}.\]
\end{theorem}

There are two parts to this theorem, one part saying that for $p < c n^{-1/m_2(H)}$ the random graph $G_{n,p}$ is highly unlikely to be $(H, r)$-Ramsey and the other saying that for $p > C n^{-1/m_2(H)}$ it is almost surely $(H,r)$-Ramsey. Following standard usage, we will refer to these two parts as the $0$-statement and the $1$-statement, respectively. 

The threshold in Theorem~\ref{thm:RR} occurs at $p^* = n^{-1/m_2(H)}$. This is the largest probability for which there is some subgraph $H'$ of $H$ such that the number of copies of $H'$ in $G_{n,p}$ is approximately the same as the number of edges. For $p$ significantly smaller than $p^*$, the number of copies of $H'$ will also be significantly smaller than the number of edges. A rather delicate argument \cite{RR93} then allows one to show that the edges of the graph may be colored so as to avoid any monochromatic copies of $H'$. For $p$ significantly larger than $p^*$, almost every edge in the random graph is contained in many copies of every subgraph of $H$. The intuition, which takes substantial effort to make rigorous \cite{RR95}, is that these overlaps are enough to force the graph to be Ramsey.

That the proof of the $0$-statement is delicate is betrayed by the omitted cases, which have smaller thresholds. For example, if a graph contains the star $K_{1, r(t-1) + 1}$, then any $r$-colouring of the edges of this graph will contain a monochromatic $K_{1,t}$. However, the threshold for the appearance of $K_{1, r(t-1) + 1}$ is lower than the threshold suggested by $m_2(K_{1,t})$. A more subtle case is when $H = P_4$, the path with $3$ edges (and $4$ vertices), and $r=2$. In this case, a cycle of length five with a pendant edge at each vertex is $(P_4, 2)$-Ramsey. While the threshold for the appearance of these graphs is at $n^{-1}$, which is the same as $n^{-1/m_2(H)}$, the threshold is coarse. This means that they start to appear with positive probability already when $p = c/n$ for any positive $c$. This implies that the $0$-statement only holds when $p = o(1/n)$.

It is worth saying a little about the proof of the $1$-statement in Theorem~\ref{thm:RR}. We will focus on the case when $H = K_3$ and $r = 2$. The key idea is to write $G_{n,p}$ as the union of two independent random graphs $G_{n, p_1}$ and $G_{n,p_2}$, chosen so that 
\[p = p_1 + p_2 - p_1 p_2 \mbox{\hspace{3mm} and \hspace{3mm}} p_2 = L p_1\]
for some large constant $L$. We first expose the smaller random graph $G_{n,p_1}$. With high probability, every  colouring of $G_{n,p_1}$ will contain many monochromatic paths of length $2$. If $p_1$ is a sufficiently large multiple of  $1/\sqrt{n}$, it is also possible to show that with high probability these monochromatic paths are well distributed. In particular, for any given colouring of $G_{n,p_1}$, there are at least $cn^3$ triangles in the underlying graph $K_n$ such that there is a path of the same colour, say red, between each pair of vertices in each triangle. 

We now expose $G_{n,p_2}$. If this graph contains any of the $cn^3$ triangles described above, we are done, since each edge of this triangle must take the colour blue. Otherwise, together with the red connecting path, we would have a red triangle. By Janson's inequality \cite{J90}, the probability that $G_{n,p_2}$ does not contain any of the $c n^3$ triangles associated to this particular colouring is at most $2^{-c' p_2 n^2}$, where $c'$ depends on $c$. However, we must remember to account for every possible colouring of $G_{n,p_1}$. To do this, we take a union bound. Indeed, since there are at most $2^{p_1 n^2}$ colourings of $G_{n,p_1}$, the probability that there exists a colouring such that $G_{n,p_2}$ does not intersect the associated set of triangles is at most $2^{p_1 n^2} 2^{-c' p_2 n^2}$. If we choose $L$ sufficiently large, this probability tends to zero, completing the proof.

This method also allowed R\"odl and Ruci\'nski to determine the threshold for van der Waerden's theorem to hold in random subsets of the integers. Van der Waerden's theorem \cite{VdW27} states that for any natural numbers $k$ and $r$ there exists $n$ such that any $r$-colouring of $[n] := \{1,2,\dots, n\}$ contains a monochromatic $k$-term arithmetic progression, that is, a monochromatic subset of the form $\{a, a+d, \dots, a+(k-1)d\}$. To state the random version of this theorem, we define $[n]_p$ to be a random subset of $[n]$ where each element is chosen independently with probability $p$. We also say that a subset $I$ of the integers is $(k, r)$-vdW if in any $r$-colouring of the points of $I$ there is a monochromatic $k$-term arithmetic progression. R\"odl and Ruci\'nski's random van der Waerden theorem \cite{RR95, RR97} is then as follows.

\begin{theorem} \label{thm:vdW}
For any positive integers $k \geq 3$ and $r \geq 2$, there exist positive constants $c$ and $C$ such that 
\[
\lim_{n \rightarrow \infty} \mathbb{P} [[n]_{p} \mbox{ is $(k,r)$-vdW}] =
\begin{cases}
0 & \text{if $p < c n^{-1/(k-1)}$}, \\
1 & \text{if $p > C n^{-1/(k-1)}$}.
\end{cases}
\]
\end{theorem}

The threshold is again a natural one, since it is the point where we expect that most vertices in $[n]_p$ will be contained in a constant number of $k$-term arithmetic progressions. We will say more about this in the next section when we discuss density theorems. 

One question left open by the work of R\"odl and Ruci\'nski was whether Theorem~\ref{thm:RR} could be extended to hypergraphs. While some partial progress was made \cite{RR98, RRS07}, the general problem remained open, not least because of the apparent need to apply a hypergraph analogue of the regularity lemma, something which has only been developed in recent years \cite{G07, NRS06, RS04, T06}. Approaches which circumvent hypergraph regularity were developed independently by Friedgut, R\"odl and Schacht \cite{FRS10} and by Gowers and the author \cite{CG14}, so that the following generalisation of Theorem~\ref{thm:RR} is now known. We write $G^{(k)}_{n,p}$ for the random $k$-uniform hypergraph on $n$ vertices, where each edge is chosen independently with probability $p$.

\begin{theorem} \label{thm:RamseyHyper}
For any $k$-uniform hypergraph $H$ and any positive integer~$r \geq 2$, there exists $C > 0$ such that 
\[
\lim_{n \rightarrow \infty} \mathbb{P} [G^{(k)}_{n,p} \mbox{ is $(H,r)$-Ramsey}] =
1  \text{ if $p > C n^{-1/m_k(H)}$},
\]
where 
\[m_k(H) = \max\left\{\frac{e(H') - 1}{v(H') - k}: H' \subseteq H \mbox{ and } v(H') \geq k+1\right\}.\]
\end{theorem}

We note that the approach in \cite{CG14} applies when $H$ is strictly $k$-balanced, that is, when $m_k(H) > m_k (H')$ for every subgraph $H'$ of $H$. However, almost all hypergraphs, including the complete hypergraph $K_t^{(k)}$, satisfy this requirement. A similar caveat applies to many of the theorems stated in this survey. We will usually make this explicit. 

The $0$-statement corresponding to Theorem~\ref{thm:RamseyHyper} was considered by Gugelmann, Person, Steger and Thomas (see \cite{G14, GPST12}). In particular, their results imply the corresponding $0$-statement for complete hypergraphs. However, there are again cases where the true threshold is smaller than $n^{-1/m_k(H)}$. Indeed, the picture seems to be more complicated than for graphs since there are examples other than the natural generalisations of paths and stars for which the $1$-statement may be improved. We refer the reader to \cite{G14} for a more complete discussion.

One may also consider the threshold for asymmetric Ramsey properties. We say that a graph $G$ is $(H_1, H_2, \dots, H_r)$-Ramsey if any colouring of the edges of $G$ with colours $1, 2, \dots, r$ contains a monochromatic copy of $H_i$ in colour $i$ for some $i \in \{1, 2, \dots, r\}$. A conjecture of Kohayakawa and Kreuter \cite{KK97}, which generalises Theorem~\ref{thm:RR}, says that if $H_1, H_2, \dots, H_r$ are graphs with $1 < m_2(H_r) \leq \cdots \leq m_2(H_1)$, then the $(H_1, H_2, \dots, H_r)$-Ramsey property has a threshold at $n^{-1/m_2(H_1, H_2)}$, where
\[m_2(H_1, H_2) = \max\left\{\frac{e(H_1')}{v(H_1') - 2 + 1/m_2(H_2)}: H_1' \subseteq H_1 \mbox{ and } v(H_1') \geq 3\right\}.\]  
Since the $0$-statement fails to hold for certain forests, this statement should be qualified further, but it seems likely to hold for most collections of graphs. 

Kohayakawa and Kreuter established the conjecture when $H_1, H_2, \dots, H_r$ are cycles. As noted in \cite{MSSS09}, the same method shows that the K\L R conjecture (which we discuss in Section~\ref{sec:regularity}) would imply the $1$-statement of the conjecture when $H_1$ is strictly $2$-balanced, that is, when $m_2(H_1) > m_2(H_1')$ for all proper subgraphs $H_1'$. Since the K\L R conjecture is now an established fact, the following theorem is known to hold (as was noted explicitly by Balogh, Morris and Samotij \cite{BMS14}).

\begin{theorem} \label{thm:asymRam}
For any graphs $H_1, H_2, \dots, H_r$ with $1 < m_2(H_r) \leq \cdots \leq m_2(H_1)$ and such that $H_1$ is strictly $2$-balanced, there exists $C > 0$ such that
\[
\lim_{n \rightarrow \infty} \mathbb{P} [G_{n,p} \mbox{ is $(H_1, H_2, \dots, H_r)$-Ramsey}] =
1  \text{ if $p > C n^{-1/m_2(H_1, H_2)}$}.
\]
\end{theorem}

A slightly weaker statement was established by Kohayakawa, Schacht and Sp\"{o}hel \cite{KSS14} without appealing to the K\L R conjecture. Their proof is much closer in spirit to R\"odl and Ruci\'nski's proof of Theorem~\ref{thm:RR}. A corresponding $0$-statement when $H_1, H_2, \dots, H_r$ are cliques was established by Marciniszyn, Skokan, Sp\"ohel and Steger \cite{MSSS09}. However, the $0$-statement remains open in general. 

The methods developed in \cite{CG14} and \cite{FRS10} also allow one to extend R\"odl and Ruci\'nski's results on random analogues of van der Waerden's theorem to a more general setting. A classical theorem of Rado \cite{R41} generalises van der Waerden's theorem by establishing necessary and sufficient conditions for a system of homogeneous linear equations 
\[\sum_{j=1}^k a_{ij} x_j = 0 \mbox{ for $1 \leq i \leq \ell$}\]
to be partition regular, that is, to be such that any finite colouring of the natural numbers contains a monochromatic solution $(x_1, x_2, \dots, x_k)$ to this system of equations. To give an example, the solutions to the system of equations $x_i + x_{i+2} = 2 x_{i+1}$ for $i = 1, 2, \dots, k-2$ are $k$-term arithmetic progressions and so van der Waerden's theorem implies that this system of equations is partition regular. An extension of Theorem~\ref{thm:vdW} was proved by R\"odl and Ruci\'nski in \cite{RR97}, but their $1$-statement only applied to density regular systems of equations (though see also \cite{GRR96}). These are systems of equations, like the system defining $k$-term arithmetic progressions, whose solutions sets are closed under translation and dilation. 

An extension of this theorem which applies to all partition regular systems of equations was proved by Friedgut, R\"odl and Schacht \cite{FRS10}. More precisely, they proved a $1$-statement, while the $0$-statement had been established earlier by R\"odl and Ruci\'nski \cite{RR97}. Since the details are somewhat technical, we refer the interested reader to \cite{FRS10} for further particulars.

We have already mentioned that the result of Frankl and R\"odl \cite{FR86} may be used to prove that there are $K_4$-free graphs which are $(K_3,2)$-Ramsey. This was originally proved by Folkman \cite{F70} using a constructive argument. More generally, he proved that for any positive integer $t$ there is a $K_{t+1}$-free graph which is $(K_t, 2)$-Ramsey. This beautiful result was subsequently extended to $r$-colourings by Ne\v set\v ril and R\"odl \cite{NR76, NR81}.

Once we know that these graphs exist, it is natural to try and estimate their size. We define the Folkman number $F(t)$ to be the smallest natural number $n$ such that there exists a $K_{t + 1}$-free graph $G$ on $n$ vertices with the property that every $2$-colouring of the edges of $G$ contains a monochromatic $K_{t}$. The upper bounds on $F(t)$ which come from the constructive proofs tend to have a dependency on $t$ which, with a conservative estimate, is at least tower-type, that is, a tower of twos of height at least $t$. On the other hand, the lower bound is essentially the same as for Ramsey's theorem, that is, $F(t) \geq 2^{c' t}$. 

Very recently, it was noted that some of the methods for proving Ramsey-type theorems in random sets yield significantly stronger bounds for Folkman numbers \cite{CG142, RRS14}. In particular, the following result was proved by R\"odl, Ruci\'nski and Schacht \cite{RRS14}. Their proof relies heavily on the hypergraph container results developed by Balogh, Morris and Samotij~\cite{BMS14} and Saxton and Thomason~\cite{ST14} and an observation of Nenadov and Steger \cite{NS14} that allows one to apply this machinery in the Ramsey setting.

\begin{theorem} \label{Folkman}
There exists a constant $c$ such that
\[F(t) \leq 2^{c t^4 \log t}.\]
\end{theorem}

This bound is tantalisingly close to the lower bound and it would be of great interest to improve it further. Since we have now brought our discussions of Ramsey-type theorems in random sets full circle, this provides a convenient departure point to move on to discussing density theorems in random sets, a topic about which much less was known before recent developments.

\section{Density theorems in random sets} \label{sec:density}

Tur\'an's theorem \cite{T41} states that the largest $K_t$-free subgraph of $K_n$ has at most $\left(1 - \frac{1}{t-1}\right) \frac{n^2}{2}$ edges. Moreover, the unique $K_t$-free subgraph achieving this maximum is the $(t-1)$-partite graph with vertex sets $V_1, V_2, \dots, V_{t-1}$, where each set is of order $\lfloor \frac{n}{t-1} \rfloor$ or $\lceil \frac{n}{t-1} \rceil$. In particular, for $t = 3$, the triangle-free subgraph of $K_n$ with the most edges is a bipartite graph with parts of order $\lfloor \frac{n}{2} \rfloor$ or $\lceil \frac{n}{2} \rceil$. A substantial generalisation of this theorem, known as the Erd\H{o}s--Stone--Simonovits theorem \cite{ESi66, ES46}, states that for any graph $H$ the largest $H$-free subgraph of $K_n$ has at most $\left(1 - \frac{1}{\chi(H) - 1} + o(1) \right) \binom{n}{2}$ edges, where $\chi(H)$ is the chromatic number of $H$.

We say that a graph $G$ is $(H, \epsilon)$-Tur\'an if every subgraph of $G$ with at least $\left(1 - \frac{1}{\chi(H)-1} + \epsilon\right) e(G)$ edges contains a copy of $H$. The original work of Frankl and R\"odl \cite{FR86} on Ramsey properties in random graphs was actually motivated by a problem of Erd\H{o}s and Ne\v set\v ril concerning an analogue of Folkman's theorem for the $(H, \epsilon)$-Tur\'an property. Specifically, they asked whether there exist $K_4$-free graphs which are $(K_3, \epsilon)$-Tur\'an and Frankl and R\"odl showed that there are. Though not stated explicitly in their paper, Frankl and R\"odl's method implies that for any $\epsilon > 0$ there exists $C > 0$ such that if $p > C/\sqrt{n}$ then
\[\lim_{n \rightarrow \infty} \mathbb{P}[G_{n,p} \mbox{ is $(K_3, \epsilon)$-Tur\'an}] =1.\]
Unlike Ramsey properties, the corresponding $0$-statement is easy to prove. Indeed, for $p$ a sufficiently small multiple of $1/\sqrt{n}$, the number of triangles in $G_{n,p}$ will be significantly smaller than the number of edges. We may therefore remove all copies of $K_3$ by deleting one edge from each copy, leaving a subgraph which is triangle-free but contains at least $(1 - \delta) e(G_{n,p})$ edges.

A similar argument provides a lower bound for all $H$. That is, if the number of copies of $H$ is significantly smaller than the number of edges, we can remove all copies of $H$ by deleting one edge from each copy. Therefore, if $p^{e(H)} n^{v(H)} \ll p n^2$, that is, $p \ll n^{-(v(H)-2)/(e(H)-1)}$, the $(H, \epsilon)$-Tur\'an property cannot hold. Since the same argument applies for any subgraph $H'$ of $H$, it is easy to see that for $p \ll n^{-1/m_2(H)}$ the random graph $G_{n,p}$ cannot be $(H, \epsilon)$-Tur\'an. Here $m_2(H)$ is defined as in Theorem~\ref{thm:RR}, that is,
\[m_2(H) = \max\left\{\frac{e(H') - 1}{v(H') - 2}: H' \subseteq H \mbox{ and } v(H') \geq 3\right\},\]
The natural conjecture that the $(H, \epsilon)$-Tur\'an property holds in random graphs with $p \gg n^{-1/m_2(H)}$ was first stated by Haxell, Kohayakawa and \L uczak \cite{HKL95, HKL96} and reiterated by Kohayakawa, \L uczak and R\"odl \cite{KLR97}.

Until recently, this conjecture was only known to hold for a small collection of graphs, including $K_3$, $K_4$ and $K_5$ \cite{FR86, KLR97, GSS04} and all cycles \cite{F94, HKL95, HKL96} (see also \cite{KRS04, SV03}). A verification of the conjecture for all graphs was completed by Schacht \cite{S14} and by Gowers and the author \cite{CG14}, although we must qualify this statement by saying that the results of \cite{CG14} apply when $H$ is strictly $2$-balanced, that is, when $m_2(H') < m_2(H)$ for all $H' \subset H$. However, the class of strictly $2$-balanced graphs includes many of the graphs one normally considers, such as cliques and cycles.

\begin{theorem} \label{thm:randomTuran}
For any graph $H$ and any $\e > 0$, there exist positive
constants $c$ and $C$ such that
\[
\lim_{n \rightarrow \infty} \mathbb{P} [G_{n,p} \mbox{ is $(H,\e)$-Tur{\'a}n}] =
\begin{cases}
0 & \text{if $p < c n^{-1/m_2(H)}$}, \\
1 & \text{if $p > C n^{-1/m_2(H)}$}.
\end{cases}
\]
\end{theorem}

As mentioned in the introduction, Schacht's proof of Theorem~\ref{thm:randomTuran} builds on R\"odl and Ruci\'nski's proof of Theorem~\ref{thm:RR}. In the last section, we gave a brief description of their method, showing how it was best to think of the random graph $G_{n,p}$ as the union of two independent random graphs $G_{n,p_1}$ and $G_{n,p_2}$. In Schacht's method, this multi-round exposure is taken further, the rough idea being to expose $G_{n,p}$ over several successive rounds and to apply a density increment argument. 

The method employed in \cite{CG14} relies upon proving a transference principle, an idea which originates in the work of Green and Tao \cite{GT08} (see also \cite{G08, RTTV08}). In the case of triangles, this transference principle says that for $p \geq C/\sqrt{n}$ any subgraph $G$ of $G_{n,p}$ may be modelled by a subgraph $K$ of the complete graph $K_n$ in such a way that the proportion of edges and triangles in $K$ is close to the proportion of edges and triangles in $G$. That is, if the sparse graph $G$ contains $c_1 p n^2$ edges and $c_2 p^3 n^3$ triangles, then the dense model $K$ will contain approximately $c_1 n^2$ edges and $c_2 n^3$ triangles.

Suppose now that we wish to prove Tur\'an's theorem for triangles relative to a random graph. Given a subgraph $G$ of $G_{n,p}$ with $\left(\frac{1}{2} + \epsilon\right) p \binom{n}{2}$ edges, we know, once our approximation is sufficiently good, that its dense model $K$ has at least $\left(\frac{1}{2} + \frac{\epsilon}{2}\right) \binom{n}{2}$ edges. A robust version of Tur\'an's theorem \cite{ESi83} then implies that $K$ contains at least $c n^3$ triangles for some $c > 0$ depending on $\epsilon$. Provided again that our approximation is sufficiently good, this implies that $G$ contains at least $\frac{c}{2} p^3 n^3$ triangles, which is even more than we required.

Though the analogue of Tur\'an's theorem for hypergraphs is rather poorly understood (see, for example, \cite{K11}), a similar strategy shows that it is still possible to transfer it to the random setting. To state the result, we need some definitions. Given a $k$-uniform hypergraph $H$, we let ex$(n, H)$ be the largest number of edges in an $H$-free subgraph of $K_n^{(k)}$ and
\[\pi_k(H) = \lim_{n \rightarrow \infty} \frac{\mbox{ex}(n,H)}{\binom{n}{k}}.\]
We then say that a $k$-uniform hypergraph $G$ is $(H, \epsilon)$-Tur\'an if every subgraph of $G$ with at least $\left(\pi_k(H) + \epsilon\right) e(G)$ edges contains a copy of $H$. Let $m_k(H)$ be defined as in the previous subsection, that is,
\[m_k(H) = \max\left\{\frac{e(H') - 1}{v(H') - k}: H' \subseteq H \mbox{ and } v(H') \geq k+1\right\}.\]
Then the analogue of Theorem~\ref{thm:randomTuran}, proved in \cite{CG14, S14}, states that the property of being $(H, \epsilon)$-Tur\'an for a $k$-uniform hypergraph $H$ has a threshold at $n^{-1/m_k(H)}$.

\begin{theorem} \label{thm:randomTuranHyper}
For any $k$-uniform hypergraph $H$ and any $\e > 0$, there exist positive
constants $c$ and $C$ such that
\[
\lim_{n \rightarrow \infty} \mathbb{P} [G^{(k)}_{n,p} \mbox{ is $(H,\e)$-Tur{\'a}n}] =
\begin{cases}
0 & \text{if $p < c n^{-1/m_k(H)}$}, \\
1 & \text{if $p > C n^{-1/m_k(H)}$}.
\end{cases}
\]
\end{theorem}

One structural counterpart to Tur\'an's theorem is the Erd\H{o}s-Simonovits stability theorem \cite{Si68}. This says that for any graph $H$ with $\chi(H) \geq 3$ and any $\epsilon > 0$, there exists $\delta > 0$ such that any $H$-free subgraph of $K_n$ with at least $\left(1 - \frac{1}{\chi(H) - 1} - \delta\right) \binom{n}{2}$ edges may be made $(\chi(H) - 1)$-partite by removing at most $\epsilon n^2$ edges. The following sparse analogue of this result was originally proved in \cite{CG14} for strictly $2$-balanced graphs. Later, Samotij \cite{Sj14} found a way to amend Schacht's method so that it applied to stability statements, extending this result to all graphs.

\begin{theorem} \label{thm:stab}
For any graph $H$ with $\chi(H) \geq 3$ and any $\epsilon > 0$, there exist positive constants $\delta$ and $C$ such that if $p \geq C n^{-1/m_2(H)}$ the random graph $G_{n,p}$ a.a.s.~has the following property. Every $H$-free subgraph of $G_{n,p}$ with at least $\left(1 - \frac{1}{\chi(H)-1} - \delta\right)p\binom{n}{2}$ edges can be made $(\chi(H)-1)$-partite by removing at most $\epsilon p n^2$ edges. 
\end{theorem}

For cliques, Tur\'an's theorem has a much more precise corresponding structural statement, saying that the largest $K_t$-free subgraph is $(t-1)$-partite. One may therefore ask when this property holds a.a.s.~in the random graph $G_{n,p}$. This question was first studied by Babai, Simonovits and Spencer \cite{BSS90} who showed that for $p > \frac{1}{2}$ the size of the maximum triangle-free subgraph is a.a.s.~the same as the size of the largest bipartite subgraph. This result was extended to the range $p > n^{-c}$ by Brightwell, Panagiotou and Steger \cite{BPS12}. Recently, DeMarco and Kahn \cite{DK14} proved the following much more precise result.

\begin{theorem} \label{thm:KDeM}
There is a positive constant $C$ such that if $p > C \sqrt{\log n/n}$ then a.a.s.~every maximum triangle-free subgraph of $G_{n,p}$ is bipartite.
\end{theorem}

The threshold here is different from the $1/\sqrt{n}$ we have come to expect. However, the result is sharp up to the constant $C$. Indeed, for $p = 0.1 \sqrt{\log n/n}$, the random graph $G_{n,p}$ will typically contain a $5$-cycle none of whose edges are contained in a triangle. In a forthcoming paper, DeMarco and Kahn \cite{DK142} prove the following extension of this result to all cliques. Once again, the extra log factors are essential.

\begin{theorem} \label{thm:KDem2}
For any natural number $t$, there exists $C > 0$ such that if 
\[p > C n^{-\frac{2}{t+1}} \log^{\frac{2}{(t+1)(t-2)}} n\] 
then a.a.s.~every maximum $K_t$-free subgraph of $G_{n,p}$ is $(t-1)$-partite. 
\end{theorem}

We note that a related question, where one wishes to determine the range of $m$ for which most $K_t$-free graphs with $n$ vertices and $m$ edges are $(t-1)$-partite, was solved recently by Balogh, Morris, Samotij and Warnke~\cite{BMSW14}. 

%The behaviour in this case also changes at around the same point as Theorem~\ref{thm:KDem2}. 

%namely, when $m \approx n^{2 - \frac{2}{t+1}} \log^{\frac{2}{(t+1)(t-2)}} n$.

The methods of \cite{CG14} and \cite{S14} also allow one to prove sparse analogues of density statements from other settings. For example, Szemer\'edi's theorem \cite{Sz75} states that for any natural number $k$ and any $\delta >0$ there exists $n_0$ such that if $n \geq n_0$ any subset of $[n]$ of density at least $\delta$ contains a $k$-term arithmetic progression. This is the density version of van der Waerden's theorem and trivially implies that theorem by taking $\delta = \frac{1}{r}$ and considering the largest colour class. This theorem and the tools arising in its many proofs \cite{F77, G01, R14} have been enormously influential in the development of modern combinatorics.

We say that a subset $I$ of the integers is $(k,\delta)$-Szemer\'edi if any subset of $I$ with at least $\delta |I|$ elements contains an arithmetic progression of length $k$. Szemer\'edi's theorem says that for $n$ sufficiently large the set $[n]$ is $(k, \delta)$-Szemer\'edi, while a striking corollary of Green and Tao's work on arithmetic progressions in the primes \cite{GT08} says that for $n$ sufficiently large the set of primes up to $n$ is $(k, \delta)$-Szemer\'edi. 

For random subsets of the integers, the $(k, \delta)$-Szemer\'edi property was first studied by Kohayakawa, \L uczak and R\"odl \cite{KLR96}, who proved that the property of being $(3, \delta)$-Szemer\'edi has a threshold at $1/\sqrt{n}$. In general, the natural conjecture is that the $(k, \delta)$-Szemer\'edi property has a threshold at $n^{-1/(k-1)}$. The lower bound is again straightforward, since for $p \ll n^{-1/(k-1)}$ the number of $k$-term arithmetic progressions is significantly smaller than the number of elements in the random set $[n]_p$, allowing us to remove one element from each arithmetic progression without significantly affecting the density. The corresponding $1$-statement was provided in \cite{CG14} and \cite{S14}.

\begin{theorem} \label{thm:Szem}
For any integer $k \geq 3$ and $\delta > 0$, there exist positive constants $c$ and $C$ such that 
\[
\lim_{n \rightarrow \infty} \mathbb{P} [[n]_{p} \mbox{ is $(k,\delta)$-Szemer\'edi}] =
\begin{cases}
0 & \text{if $p < c n^{-1/(k-1)}$}, \\
1 & \text{if $p > C n^{-1/(k-1)}$}.
\end{cases}
\]
\end{theorem}

A particularly satisfying approach to density theorems in random sets is provided by the recent hypergraph containers method of Balogh, Morris and Samotij \cite{BMS14} and Saxton and Thomason \cite{ST14}, the only probabilistic input being Chernoff's inequality and the union bound. In the context of Szemer\'edi's theorem, one of the main corollaries of this method is the following theorem.

\begin{theorem} \label{thm:HCSzem}
For any integer $k \geq 3$ and any $\epsilon > 0$, there exists $C > 0$ such that if $m \geq C n^{1 - 1/(k-1)}$, then there are at most $\binom{\epsilon n}{m}$ subsets of $\{1, 2, \dots, n\}$ of order $m$ which contain no $k$-term arithmetic progression. 
\end{theorem}

Given this statement, which is completely deterministic, it is straightforward to derive the $1$-statement in Theorem~\ref{thm:Szem}, so much so that we may now give the entire calculation. For brevity, we write $(k, \delta)$-Sz rather than $(k, \delta)$-Szemer\'edi and $\mathcal{I}_k(n, \delta pn/2)$ for the collection of subsets of $\{1, 2, \dots, n\}$ of order $\delta pn/2$ which contain no $k$-term arithmetic progression. We have 
\begin{align*}
\mathbb{P}[[n]_p \mbox{ is not $(k, \delta)$-Sz}] & \leq \mathbb{P}[|[n]_p| < pn/2] + \mathbb{P}[|[n]_p| \geq pn/2 \mbox { and $[n]_p$ is not $(k, \delta)$-Sz}]\\
& \leq \exp(-\Omega(pn)) + \mathbb{P}[[n]_p \supseteq I \mbox{ for some } I \in \mathcal{I}_k(n, \delta pn/2)]\\
& \leq \exp(-\Omega(pn)) + \binom{\epsilon n}{\delta pn/2} p^{\delta pn/2}\\
& \leq \exp(-\Omega(pn)) + \left(\frac{2e \epsilon pn}{\delta pn}\right)^{\delta pn/2}\\
& = \exp(-\Omega(pn)),
\end{align*}
provided $\epsilon < \delta/2 e$. 

Deriving Theorem~\ref{thm:randomTuran} from the results of \cite{BMS14} and \cite{ST14} involves a little more work. To describe the idea, we focus on the case where $H = K_3$. We begin by considering the $3$-uniform hypergraph $\mathcal{G}$ whose vertex set is the collection of edges in $K_n$ and whose edge set is the collection of triangles in $K_n$. Tur\'an's theorem for triangles may then be restated as saying that this $3$-uniform hypergraph has no independent set of order greater than $\left(\frac{1}{2} + o(1)\right) |V(\mathcal{G})|$. We would now like to show that if $p \geq C/\sqrt{n}$ then the random set $V(\mathcal{G})_p$ formed by choosing each element of $V(\mathcal{G})$ independently with probability $p$ contains no independent set of order greater than $\left(\frac{1}{2} + \epsilon\right) p |V(\mathcal{G})|$. 

One approach would be to use the union bound and Chernoff's inequality to show that with high probability the intersection of the random set with each independent set is as required. An argument of this variety worked in the proof of Theorem~\ref{thm:Szem} above, but usually there are far too many independent sets for this approach to be feasible. The main results in both \cite{BMS14} and \cite{ST14} circumvent this difficulty by showing that there is a substantially smaller collection of almost independent sets which contain all independent sets. Since these sets are almost independent, we know, by the robust version of Tur\'an's theorem, that they must also have size at most $\left(\frac{1}{2} + \frac{\epsilon}{2}\right)|V(\mathcal{G})|$, say. Applying the union bound over this smaller set then allows us to derive the result.

\section{Regularity in random graphs} \label{sec:regularity}

Szemer\'edi's regularity lemma~\cite{Sz78} is one of the cornerstones of modern graph theory (see~\cite{KSSS02, RS10}). Roughly speaking, it says that the vertex set of every graph $G$ may be divided into a bounded number of parts in such a way that most of the induced bipartite graphs between different parts are pseudorandom. To be more precise, we need some definitions. 

We say that a bipartite graph between sets $U$ and $V$ is $\e$-regular if, for every $U' \subseteq U$ and $V' \subseteq V$ with $|U'| \geq \e |U|$ and $|V'| \geq \e |V|$, the density $d(U', V')$ of edges between $U'$ and $V'$ satisfies
\[|d(U', V') - d(U, V)| \leq \e.\]
A partition of the vertex set of a graph into $t$ pieces $V_1, \dots, V_t$ is an equipartition if, for every $1 \leq i, j \leq t$, we have $||V_i| - |V_j|| \leq 1$. Finally, a partition is $\e$-regular if it is an equipartition and, for all but at most $\e t^2$ pairs $(V_i, V_j)$, the induced graph between $V_i$ and $V_j$ is $\e$-regular. Szemer\'edi's regularity lemma can now be  stated as follows.

\begin{theorem} \label{thm:reglemma}
For any $\e > 0$, there exists an integer $T$ such that every graph $G$ admits an $\e$-regular partition $V_1, \dots, V_t$ of its vertex set into $t \leq T$ pieces. 
\end{theorem}

For sparse graphs -- that is, graphs with $n$ vertices and $o(n^2)$ edges -- the regularity lemma stated above is vacuous, since every equipartition into a bounded number of parts is $\e$-regular for $n$ sufficiently large. However, as observed independently by Kohayakawa \cite{K97} and R\"odl, there is a meaningful analogue of the regularity lemma for sparse graphs, provided one is willing to restrict consideration to a well-behaved class of graphs. 

To make this more precise, we say that a bipartite graph between sets $U$ and $V$ is $(\e, p)$-regular if, for every $U' \subseteq U$ and $V' \subseteq V$ with $|U'| \geq \e |U|$ and $|V'| \geq \e |V|$, the density $d(U', V')$ of edges between $U'$ and $V'$ satisfies
\[|d(U', V') - d(U, V)| \leq \e p.\]
That is, we alter the definition of regularity so that it is relative to a particular density $p$, usually chosen to be comparable to the total density between $U$ and $V$. A partition of the vertex set of a graph into $t$ pieces $V_1, \dots, V_t$ is then said to be $(\e, p)$-regular if it is an equipartition and, for all but at most $\e t^2$ pairs $(V_i, V_j)$, the induced graph between $V_i$ and $V_j$ is $(\e, p)$-regular. 

The class of graphs to which the Kohayakawa-R\"odl regularity lemma applies are the so-called upper-uniform graphs~\cite{KoRo03}. Suppose that $0 < \eta \leq 1$, $D > 1$ and $0 < p \leq 1$ are given. We will say that a graph $G$ is $(\eta, p, D)$-upper-uniform if for all disjoint subsets $U_1$ and $U_2$ with $|U_1|, |U_2| \geq \eta |V(G)|$, the density of edges between $U_1$ and $U_2$ satisfies $d(U_1, U_2) \leq D p$. This condition is satisfied for many natural classes of graphs, including all subgraphs of random graphs of density $p$. The sparse regularity lemma of Kohayakawa and R\"odl is now as follows.

\begin{theorem} \label{thm:sparsereg}
For any $\e > 0$ and $D > 1$, there exists $\eta > 0$ and an integer $T$ such that for every $p\in[0,1]$, every graph $G$ that is $(\eta, p, D)$-upper-uniform admits an $(\e, p)$-regular partition $V_1, \dots, V_t$ of its vertex set into $t \leq T$ pieces. 
\end{theorem}

A recent variant of this lemma, due to Scott \cite{Sc11}, requires no upper-uniformity assumption on $G$, although it is often useful to impose such a constraint in practice. Since the two statements are interchangeable when one is dealing with a subgraph of the random graph, we have chosen to describe the original version.

In applications, the regularity method is usually applied in combination with a counting lemma. Roughly speaking, a counting lemma says that if we start with an arbitrary graph $H$ and replace its vertices by large independent sets and its edges by $\e$-regular bipartite graphs with non-negligible density, then this blow-up will contain roughly the expected number of copies of $H$. To state this result formally, we again need some definitions.

Given a graph $H$ with vertex set $\{1, 2, \dots, k\}$ and a collection of disjoint vertex sets $V_1, V_2, \dots, V_k$ in a graph $G$, we say that a $k$-tuple $(v_1,v_2,\dots,v_k)$ is a canonical copy of $H$ in $G$ if $v_i\in V_i$ for every $i\in V(H)$ and $v_iv_j\in E(G)$ for every $ij\in E(H)$. We write $G(H)$ for the number of canonical copies of $H$ in $G$. The counting lemma may now be stated as follows.

\begin{lemma}
  \label{lemma:H-count}
  For any graph $H$ with vertex set $\{1,2,\dots,k\}$ and any $\delta>0$, there exists a positive constant $\e$ and an integer $n_0$ such that the following holds. Let $n\geq n_0$ and let $G$ be a graph whose vertex set is a disjoint union $V_1 \cup V_2 \cup \dots \cup V_k$ of sets of size $n$. Assume that for each $ij \in E(H)$, the bipartite subgraph of $G$ induced by $V_i$ and $V_j$ is $\e$-regular and has density $d_{ij}$. Then  
\[G(H) = \left(\prod_{ij\in E(H)}d_{ij}\pm\delta\right) n^k.\]
\end{lemma}

When combined with the regularity lemma, this result allows one to prove a number of well-known theorems in extremal graph theory, including the Erd\H{o}s--Stone--Simonovits theorem, its stability version and the graph removal lemma. In order to extend these results to sparse graphs, one plausible approach, championed by Kohayakawa, \L uczak and R\"odl \cite{KLR97}, would be to extend Lemma~\ref{lemma:H-count} to sparse graphs. For example, it would be ideal if we could replace the densities $d_{ij}$ with $d_{ij} p$, the $\epsilon$-regularity condition with an $(\epsilon, p)$-regularity condition and the conclusion with 
\[G(H) =  \left(\prod_{ij\in E(H)}d_{ij}\pm\delta\right) p^{e(H)} n^k.\]
We will initially aim for less, only asking to embed a single canonical copy of $H$. Unfortunately, for reasons with which we are now familiar, we cannot hope that such a statement holds for small $p$. Indeed, if $p \ll n^{-1/m_2(H)}$, there is a subgraph $H'$ of $H$ for which $p^{e(H')} n^{v(H')} \ll p n^2$. We may therefore remove all copies of $H'$, and hence $H$, from $G_{n,p}$ while deleting only a small fraction of the edges. The resulting graph is both $(\epsilon, p)$-regular, for some small $\epsilon$, and $H$-free. 

Frustratingly, this embedding lemma also fails for larger values of $p$. To see this, take a counterexample of the kind just described but with the sets $V_i$ of order $r$ for some $r$ that is much smaller than $n$. Now replace each vertex of this small graph by an independent set with $n/r$ vertices and each edge with a complete bipartite graph. This yields a graph with $n$ vertices in each $V_i$. It is easy to see that the counterexample survives this blowing-up process, implying that the sought-after sparse embedding lemma is false whenever $p=o(1)$ (see \cite{GS05, KR03}).

However, these counterexamples have a very special structure, an observation that led Kohayakawa, \L uczak and R\"odl to conjecture that they might be rare. Roughly speaking, their conjecture, known as the K\L R conjecture, stated that if $p \gg n^{-1/m_2(H)}$, then the number of counterexamples to the embedding lemma is so small that $G_{n,p}$ should not typically contain any such counterexample as a subgraph. Before stating the conjecture (or theorem as it is now), we introduce some notation. 

As above, let $H$ be a graph with vertex set $\{1,2,\dots,k\}$. We denote by $\G(H, n, m, p, \e)$ the collection of all graphs $G$ obtained in the following way. The vertex set of $G$ is a disjoint union $V_1 \cup V_2 \cup \dots \cup V_k$ of sets of size $n$.  For each edge $ij \in E(H)$, we add an $(\e,p)$-regular bipartite graph with $m$ edges between the pair $(V_i, V_j)$. These are the only edges of $G$. We also write $\G^*(H,n,m,p,\e)$ for the set of all $G\in\G(H, n, m, p, \e)$ that do not contain a canonical copy of $H$.

Since the sparse regularity lemma could yield graphs with different densities between the various pairs of vertex sets, it may seem surprising that we are restricting attention to graphs where all the densities are equal. However, it is sufficient to consider just this case. In fact, the K\L R conjecture, which we now state, is more specific still, since it also takes $p = m/n^2$. Again, it turns out that from this case one can deduce any other cases that may be needed.

\begin{theorem} \label{thm:KLR}
Let $H$ be a fixed graph and let $\beta>0$. Then there exist positive constants $C$ and $\e$ such that 
\[|\G^*(H,n,m,m/n^2,\e)| \leq \beta^m \binom{n^2}{m}^{e(H)}\]
for every $m \geq C n^{2- 1/m_2(H)}$.
\end{theorem}

The K\L R conjecture has attracted considerable attention over the past two decades and was resolved for a number of special cases. The cases $H = K_3$, $K_4$ and $K_5$ were solved in \cite{KLR96}, \cite{GPSST07}, and \cite{GSS04}, respectively. For cycles, the conjecture was proved in~\cite{B02, GKRS07} (see also~\cite{KK97} for a slightly weaker version). Related results were also given in \cite{GMS07} and \cite{KRS04}. We state it as a theorem because it has now been proved in full generality by Balogh, Morris and Samotij~\cite{BMS14} and by Saxton and Thomason~\cite{ST14}.

Many of the results discussed in this survey, including Theorems~\ref{thm:randomTuran} and \ref{thm:stab}, follow easily from the K\L R conjecture. Indeed, these applications were the original motivation for the conjecture. However, there are situations where an embedding result is not enough: rather than just a single copy of $H$, one needs to know that there are many copies. That is, one needs something more like a full counting lemma. Such a counting lemma was provided in a paper of Gowers, Samotij, Schacht and the author \cite{CGSS14}, the main result of which is the following. We allow for different densities between parts by replacing $m$ with a vector $\mathbf{m} = (m_{ij})_{ij \in E(H)}$.

\begin{theorem}
  \label{thm:KLRCount}
  For any graph $H$ and any $\delta, d > 0$, there exist positive constants $\e$ and $\xi$ with the following property. For any $\eta > 0$, there is $C > 0$ such that if $p \geq C N^{-1/m_2(H)}$ then a.a.s.\ the following holds in $G_{N,p}$:
  \begin{enumerate}
  \item[(i)]
    For any $n \geq \eta N$, $\mathbf{m}$ with $m_{ij} \geq d p n^2$ for all $ij \in E(H)$ and any subgraph $G$ of $G_{N,p}$ in $\G(H,n,\mathbf{m}, p, \e)$,
    \begin{equation*} 
      G(H) \geq \xi \left( \prod_{ij \in E(H)}\frac{m_{ij}}{n^2} \right) n^{v(H)}.
    \end{equation*}
    
  \item[(ii)]
    Moreover, if $H$ is strictly $2$-balanced, then 
    \begin{equation*}
      G(H) = (1 \pm \delta) \left( \prod_{ij \in E(H)}\frac{m_{ij}}{n^2} \right) n^{v(H)}.
    \end{equation*}
        
  \end{enumerate}
\end{theorem}

We note that Theorem~\ref{thm:KLRCount}(i) follows from Samotij's adaptation \cite{Sj14} of Schacht's method \cite{S14} (and may also be derived from the work of Saxton and Thomason \cite{ST14}), while Theorem~\ref{thm:KLRCount}(ii) follows from the work of Gowers and the author \cite{CG14}. Though stronger than Theorem~\ref{thm:KLR} in some obvious ways, it is worth noting that Theorem~\ref{thm:KLRCount} does not return the estimate for the number of counterexamples provided by that theorem. This estimate is important for some applications, Theorem~\ref{thm:asymRam} being a notable example.

One sample application where we need a counting result rather than an embedding result is for proving a random analogue of the graph removal lemma. This theorem, usually attributed to Ruzsa and Szemer\'edi~\cite{RS78} (though see also~\cite{ADLRY94, EFR86, Fu95}), is as follows: for any $\delta > 0$, there exists $\epsilon > 0$ such that if $G$ is a graph on $n$ vertices containing at most $\epsilon n^{v(H)}$ copies of $H$, then $G$ may be made $H$-free by deleting at most $\delta n^2$ edges. Though simple in appearance, this result is surprisingly difficult to prove (see, for example, \cite{CF13, F11}). It also has some striking consequences, including the $k = 3$ case of Szemer\'edi's theorem, originally due to Roth~\cite{R53}. A sparse random version of the graph removal lemma was conjectured by \L uczak in \cite{L06} and proved, for strictly $2$-balanced $H$, in~\cite{CG14}. The following statement, which applies for all $H$, may be found in \cite{CGSS14}.

\begin{theorem} \label{thm:removal-Gnp}
For any graph $H$ and any $\delta > 0$, there exist positive constants $\e$ and $C$ such that if $p \geq Cn^{-1/m_2(H)}$ then the following holds a.a.s.\ in $G_{n,p}$. Every subgraph of $G_{n,p}$ which contains at most $\e p^{e(H)} n^{v(H)}$ copies of $H$ may be made $H$-free by removing at most $\delta p n^2$ edges.
\end{theorem}

Note that if $p \le c n^{-1/m_2(H)}$, for $c$ sufficiently small, this statement is trivially true. Indeed, in this range, there exists a subgraph $H'$ of $H$ such that the number of copies of $H'$ in $G_{n,p}$ is smaller than $\delta p n^2$, so we can simply remove one edge from each copy of $H'$. One might then conjecture, as \L uczak did, that Theorem \ref{thm:removal-Gnp} holds for all values of $p$. For $2$-balanced graphs, those with $m_2(H') \leq m_2(H)$ for all $H' \subset H$, we may verify this conjecture by taking $\epsilon$ to be sufficiently small in terms of $C$, $\delta$, and $H$. For $p \leq C n^{-1/m_2(H)}$ and $\epsilon < \delta C^{-e(H)}$, the number of copies of $H$ is at most $\epsilon p^{e(H)} n^{v(H)} \leq \epsilon C^{e(H)} p n^2 < \delta p n^2$. Deleting an edge from each copy yields the result.

%Note analogues for hypergraphs hold?

\section{Further directions}

\subsection{Sharp thresholds for Ramsey properties}

A graph property $\mathcal{P}$ is said to be monotone if it is closed under the addition of edges, that is, $G \in \mathcal{P}$ and $G \subset G'$ implies that $G' \in \mathcal{P}$. A result of Bollob\'as and Thomason \cite{BT87} shows that any monotone property has a threshold. For example, since Ramsey properties are clearly monotone, this immediately implies that the $(H, r)$-Ramsey property and the $(k,r)$-vdW property, both defined in Section~\ref{sec:Ramsey}, have thresholds.

Once we have proved that a given property has a threshold, it is often interesting to study this threshold more closely. We say that $\mathcal{P}$ has a sharp threshold at $p^* := p^*(n)$ if, for any $\epsilon > 0$,
\[
\lim_{n \rightarrow \infty} \mathbb{P} [G_{n,p} \mbox{ is in $\mathcal{P}$}] =
\begin{cases}
0 & \text{if $p < (1 - \epsilon)p^*$}, \\
1 & \text{if $p > (1 + \epsilon)p^*$}.
\end{cases} 
\]
For example, the properties of being connected and having a Hamiltonian cycle have sharp thresholds, while the property of containing a particular graph $H$ has a non-sharp or coarse threshold.

A seminal result of Friedgut \cite{F99} gives a criterion for assessing whether a monotone property has a sharp threshold or not. Roughly speaking, this criterion says that if the property is globally determined the threshold is sharp, while if it is locally determined it is not. This fits in well with the examples given above, since connectedness and Hamiltonicity are clearly global properties, while having a single copy of a particular $H$ is decidedly local.

The question of whether Ramsey properties have sharp thresholds was first studied by Friedgut and Krivelevich \cite{FK00}. They proved, amongst other things, that the $(H, r)$-Ramsey property is sharp when $H$ is any tree other than a star or a path of length three. However, the first substantial breakthrough was made by Friedgut, R\"odl, Ruci\'nski and Tetali \cite{FRRT06}, who proved that the $(K_3, 2)$-Ramsey property has a sharp threshold. Their result may be stated as follows.

\begin{theorem} \label{thm:FRRT}
There exists a bounded function $\hat{c} := \hat{c}(n)$ such that for any $\epsilon > 0$,
\[
\lim_{n \rightarrow \infty} \mathbb{P} [G_{n,p} \mbox{ is $(K_3,2)$-Ramsey}] =
\begin{cases}
0 & \text{if $p < (1 - \epsilon)\hat{c}/\sqrt{n}$}, \\
1 & \text{if $p > (1 + \epsilon)\hat{c}/\sqrt{n}$}.
\end{cases}
\]
\end{theorem}

A close look at this result reveals an unusual feature: though we know that the threshold is sharp, we do not know exactly where it lies. In principle, the function $\hat{c}(n)$ could depend on $n$ and wander up and down between constants $c$ and $C$. However, we expect that the true behaviour should be that it tends towards a constant. It would be very interesting to prove that this is the case. It would also be of great interest to extend Theorem~\ref{thm:FRRT} to other graphs and a higher number of colours. 

More recently, Friedgut, H\`an, Person and Schacht \cite{FHPS14} proved that there is a sharp threshold for the appearance of $k$-term arithmetic progressions in every $2$-colouring of $[n]_p$. That is, they showed that the $(k,2)$-vdW property has a sharp threshold. Their proof relies in a fundamental way on the hypergraph containers results discussed throughout this survey.

\begin{theorem}
For every integer $k \geq 3$, there exists a bounded function $\hat{c}_k := \hat{c}_k(n)$ such that for any $\epsilon > 0$,
\[
\lim_{n \rightarrow \infty} \mathbb{P} [G_{n,p} \mbox{ is $(k,2)$-vdW}] =
\begin{cases}
0 & \text{if $p < (1 - \epsilon)\hat{c}_k n^{-1/(k-1)}$}, \\
1 & \text{if $p > (1 + \epsilon)\hat{c}_k n^{-1/(k-1)}$}.
\end{cases}
\]
\end{theorem}

It would again be interesting to determine the asymptotic behaviour of $\hat{c}_k(n)$ or to extend this result to a higher number of colours.

\subsection{Large subgraph theorems in random graphs}

One of the most active areas of research in extremal combinatorics is in finding conditions under which a graph contains certain large or even spanning sparse subgraphs (see, for example, \cite{KO09}). It is therefore natural to ask whether these results also have random analogues. 

One of the standard examples in this area is Dirac's theorem \cite{D52}, which says that if a graph on $n$ vertices has minimum degree at least $n/2$ then it contains a Hamiltonian cycle, that is, a cycle which meets every vertex. The study of random analogues of Dirac's theorem was initiated by Sudakov and Vu~\cite{SV08} and the state of the art is now the following result of Lee and Sudakov~\cite{LS12}. 

\begin{theorem} \label{thm:Dirac}
For any $\epsilon > 0$, there exists $C > 0$ such that if $p \geq C\frac{\log n}{n}$ then a.a.s.~every subgraph of $G_{n,p}$ with minimum degree at least $\left(\frac{1}{2} + \epsilon \right) pn$ contains a Hamiltonian cycle.
\end{theorem}

There has also been considerable work on studying random analogues of the bandwidth theorem of B\"ottcher, Schacht and Taraz \cite{BST09}. The bandwidth of a graph $G$ is the smallest $b$ for which there is an ordering $v_1, v_2, \dots, v_n$ of the vertices of $G$ such that $|i - j| \leq b$ for all edges $v_i v_j$. The theorem then states that for any positive integers $r$ and $\Delta$ and any $\gamma > 0$, there exists an integer $n_0$ and $\beta > 0$ such that if $n \geq n_0$ and $H$ is an $n$-vertex graph with chromatic number $r$, maximum degree $\Delta$ and bandwidth at most $\beta n$, then any graph on $n$ vertices with minimum degree at least $\left(1 - \frac{1}{r} + \gamma\right) n$ contains a copy of $H$. 

For the $r = 2$ case, that is, for bipartite $H$, the following random analogue of this theorem was proved by B\"ottcher, Kohayakawa and Taraz \cite{BKT13}. 

\begin{theorem} \label{thm:bandwidth}
For any integer $\Delta \geq 2$ and any $\eta, \gamma > 0$, there exist positive constants $\beta$ and $C$ such that if $p \geq C (\log n/n)^{1/\Delta}$ the random graph $G_{n,p}$ a.a.s.~has the following property. Any subgraph of $G_{n,p}$ with minimum degree at least $\left(\frac{1}{2} + \gamma\right) p n$ contains any bipartite graph on at most $(1 - \eta)n$ vertices with maximum degree $\Delta$ and bandwidth at most $\beta n$.
\end{theorem}

Related results were also proved by Huang, Lee and Sudakov \cite{HLS12}. In particular, they showed that if $H$ is an $r$-partite graph on $n$ vertices such that every vertex is contained in a triangle, then there exist subgraphs of the random graph $G_{n,p}$ with minimum degree at least $\left(1 - \frac{1}{r} + \gamma\right) pn$ such that at least $cp^{-2}$ vertices are not contained in a copy of $H$. That is, we cannot hope to cover all vertices when considering random analogues of the bandwidth theorem. However, as suggested by results in \cite{BLS12} and \cite{HLS12}, it may still be possible to embed graphs with as many as $n - C p^{-2}$ vertices. 

A celebrated result of Chv\'atal, R\"odl, Szemer\'edi and Trotter \cite{CRST83} (see also \cite{CFS12, GRR00}) states that for any positive integers $\Delta$ and $r$, there exists $C > 0$ such that if $H$ is any graph with $n$ vertices and maximum degree $\Delta$, then $R(H; r) \leq C n$. That is, the Ramsey number of bounded degree graphs grows linearly in the number of vertices. However, one can do even better. 

Given a graph $H$ and a natural number $r$, we define the size-Ramsey number $\hat{R}(H; r)$ to be the smallest number of edges in an $(H,r)$-Ramsey graph. So we are now interested in minimising the number of edges rather than the number of vertices. A striking result of Beck \cite{B83} says that $\hat{R}(P_n;r) \leq C n$ for some $C$ depending only on $r$. Using random graphs,  Kohayakawa, R\"odl, Schacht and Szemer\'edi~\cite{KRSS11} recently proved that if $H$ is any graph with $n$ vertices and maximum degree $\Delta$, then $\hat{R}(H;r) \leq n^{2 - \frac{1}{\Delta} + o(1)}$. That is, the size-Ramsey number of bounded degree graphs is subquadratic in the number of vertices. Precisely stated, their main result is the following.

\begin{theorem} \label{thm:SizeRamsey}
For any integers $\Delta \geq 2$ and $r \geq 2$, there exists $C > 0$ such that if $p \geq C (\log N/N)^{1/\Delta}$ the random graph $G_{N,p}$ with $N = C n$ a.a.s.~has the following property. Any $r$-colouring of the edges of $G_{n,p}$ contains a colour class which contains every graph on $n$ vertices with maximum degree $\Delta$. 
\end{theorem}

In a forthcoming paper, Allen, B\"ottcher, H\`an, Kohayakawa and Person \cite{ABHKP14} prove a sparse random version of the blow-up lemma. For dense graphs, this result, proved by Koml\'os, S\'ark\"ozy and Szemer\'edi \cite{KSS97}, is a standard tool for embedding spanning subgraphs. Its sparse counterpart should allow one to reprove many of the results mentioned in this section in a unified way.

\subsection{Combinatorial theorems relative to a pseudorandom set}

While this survey has focused on combinatorial theorems relative to random sets, analogous questions may also be asked for pseudorandom sets. Much of the work in this direction has focused on the combinatorial properties of the class of $(p, \beta)$-jumbled graphs. These graphs, introduced by Thomason \cite{T87, T872}, have the property that if $X$ and $Y$ are vertex subsets, then
\[|e(X,Y) - p|X||Y|| \leq \beta \sqrt{|X||Y|}.\]
As one would expect of a pseudorandom property, the random graph $G_{n,p}$ is itself $(p,\beta)$-jumbled. In this case, with high probability, we may take $\beta$ to be  $O(\sqrt{pn})$. This is essentially optimal, that is, there are no $(p,\beta)$-jumbled graphs with $\beta = o(\sqrt{pn})$. An explicit example of a jumbled graph is the Paley graph. This is the graph with vertex set $\mathbb{Z}_p$, where $p$ is a prime of the form $4k+1$, and edge set given by joining $x$ and $y$ if and only if their difference is a quadratic residue. This graph is again optimally jumbled with $p = \frac{1}{2}$ and $\beta = O(\sqrt{n})$. For many more examples, we refer the reader to the survey \cite{KS06}. 

For $(p, \beta)$-jumbled graphs, one is usually interested in questions of the following form: given a graph property $\mathcal{P}$, an integer $n$ and a density $p$, for what values of $\beta$ is it the case that a $(p, \beta)$-jumbled graph on $n$ vertices satisfies $\mathcal{P}$? To give an example, for any integer $t \geq 3$, there exists $c > 0$ such that if $\beta \leq c p^{t-1} n$ then any $(p, \beta)$-jumbled graph on $n$ vertices contains a copy of $K_t$. For $t = 3$, this condition is known to be tight, as shown by an example of Alon \cite{Al94}. 

Very recently, a general method for transferring combinatorial theorems to pseudorandom graphs was found by Fox, Zhao and the author \cite{CFZ14}. Though we will not attempt an exhaustive survey, the following sample result is representative. 

\begin{theorem} \label{thm:PseudoRemoval}
For any integer $t$ and any $\epsilon > 0$, there exist positive constants $\delta$ and $c$ such that if $\beta \leq c p^t n$ then any $(p, \beta)$-jumbled graph $G$ on $n$ vertices has the following property. Any subgraph of $G$ containing at most $\delta p^{\binom{t}{2}} n^t$ copies of $K_t$ may be made $K_t$-free by deleting at most $\epsilon p n^2$ edges. 
\end{theorem}

That is, we have an extension of the removal lemma to subgraphs of pseudorandom graphs. Although we have only stated this result for cliques, there is also a more general statement that applies to all graphs. Moreover, with similar conditions on $\beta$, it is possible to prove analogues of many different combinatorial statements. For example, the $(K_t,r)$-Ramsey property and $(K_t, \epsilon)$-Tur\'an property both hold in pseudorandom graphs with $\beta \leq c p^t n$.

Unfortunately, there is still a gap in these results, even for triangles. For $t = 3$, Theorem~\ref{thm:PseudoRemoval} (which in this case was first proved by Kohayakawa, R\"odl, Schacht and Skokan \cite{KRSS10}) says that if $\beta \leq c p^3 n$ then the triangle removal lemma holds for subgraphs of a $(p, \beta)$-jumbled graph on $n$ vertices. However, it may well be the case that $\beta \leq c p^2 n$ is sufficient. If true, Alon's example would imply that such a result was optimal.

The method of \cite{CFZ14} was extended to hypergraphs in \cite{CFZ142}, under a different type of pseudorandomness hypothesis (though see also \cite{CFZ142sup}). This result was then used to prove a pseudorandom analogue of Szemer\'edi's theorem. Such a result was a key ingredient in Green and Tao's proof that the primes contain arbitrarily long arithmetic progressions. Their original result states that if a subset of the integers satisfies two pseudorandomness conditions, the linear forms condition and the correlation condition, then it is $(k, \delta)$-Szemer\'edi. Our results allow one to remove the correlation condition from this statement. Due to space constraints, we are unable to say more here. However, we refer the reader to \cite{CFZ143} for further details.

\vspace{3mm} \noindent
{\bf Acknowledgements.} I would like to thank Jacob Fox, Tim Gowers, Rob Morris, Wojtek Samotij, Mathias Schacht, Benny Sudakov and Yufei Zhao for helpful discussions on the topics in this paper. Particular thanks are due to Wojtek Samotij for a number of detailed comments on the manuscript.

\end{document}